\documentclass{article}
\usepackage{amsmath,amssymb,graphicx,tikz,hyperref}
\usetikzlibrary{arrows.meta,calc,angles,positioning}
\usetikzlibrary{decorations.pathmorphing}

\usetikzlibrary{patterns}

\usepackage{amsthm}
 
\newtheorem*{prop}{Proposition}
\begin{document}
%{\LARGE NOTE}\smallskip
%\hrule\bigskip

{\LARGE\bf An alternative quadratic formula}\smallskip

{\small Norbert Hungerb\"uhler%\\
%Department of Mathematics\\
%ETH Z\"urich\\
%8092 Z\"urich, Switzerland\\
%norbert.hungerbuehler@math.ethz.ch
}

\section{A short history of the quadratic formula}
The commonly used formula to solve a quadratic equation
\begin{equation}\label{abc}
ax^2+bx+c=0,\quad a,b,c\in\mathbb C, a\neq 0,
\end{equation}
is the {\em quadratic formula}
\begin{equation}\label{solabc}
x_{1,2}=\frac{-b\pm\sqrt{b^2-4ac}}{2a}.
\end{equation}
In German speaking countries this formula is called  
{\em Midnight Formula\/} or 
{\em Moonshine Formula}\footnote{{\em Mitternachtsformel\/} or {\em Mondscheinformel}}
since a student is supposed to know it by heart even
when torn from sleep at midnight. From a didactical
point of view, this paradigm is certainly disputable since 
knowing something by heart is the enemy of understanding
the subject. This is why Hermann Weyl used to urge his students
to ``un-learn'' this formula.

The formula appeared for the first time in Descartes' {\em Discours de la m\'ethode\/} \cite[p.~303]{descartes}
where he coined one of the most enduring conventions 
of mathematics, namely the use of $a, b, c,\ldots$ as constants and $x, y, z$ as 
variables (see~\cite{notices}).

The prehistory of the quadratic formula goes back to
Babylonian mathematics: One of the oldest traces
is the cuneiform tablet {\em BM 13901}, stored in the British Museum,
which states around 80 mathematical problems. It dates
back to Hammurabi's reign in the beginning of the 18th century BC.
The first problem asks to find the side length of a square (see~\cite{hoyrup} for a detailed discussion):
\begin{center}
%\parbox{.62\hsize}
{\em The surface plus the length of the square is $3/4$.}%\quad\parbox{.35\hsize}{In modern terms $x^2+x=45$}
\end{center}
In modern terms this corresponds to the equation $x^2+x=3/4$. The later Greek
mathematics would have never accepted this mixing of dimensions. However,
the proposed Babylonian solution may be translated as follows:
\begin{center}
\parbox{.92\hsize}{\em Take 1, the unity. Brake it up in $1/2$ and $1/2$.
Square $1/2$ to get $1/4$. 
Add $1/4$ to $3/4$ and obtain $1$.
The root of this is $1$. Subtract of it the $1/2$ which you squared before.
You get $1/2$ which is the side of the square.
}
\end{center}
These numerical steps are expressed in the formula
$$
\sqrt{(1/2)^2+3/4}-1/2=1/2
$$
which indeed corresponds to the positive solution of the quadratic equation (see also~(\ref{solnormal})).

In Greek mathematics we find geometric methods to solve quadratic equations
in Euclid's works: In~\cite{greek}
and~\cite{vdw} the Theorems 5 \& 6 in Book II and
the Theorems 28 \& 29 in Book VI of the {\em Elements}, 
as well as the Theorems 84 \& 85 in the {\em Data} are interpreted as
tools to solve quadratic equations (see~\cite{herz} for a critical analysis of this view).
However, it is  conceivable that Euclid would have rather used the
Theorems 35 and 36 in Book III,
since they allow a much more direct link between the quadratic equation
and its geometric solution: Suppose we want to solve the equation
$x^2-ux+vw=0$, where $u,v,w>0, u^2-4vw\ge 0$, are given as lengths of segments.
Then, the solutions $x_1,x_2$ satisfy $x_1+x_2=u$ and $x_1x_2=vw$.
Hence the following application of the intersecting chords theorem yields the solution on the left:
\begin{center}
\begin{tikzpicture}[line cap=round,>=latex,line join=round,x=2.5,y=2.5]
\clip(-40.32267637304512,-20.092518246595212) rectangle (24.328624984866302,38.38604579573683);
\draw[shorten <=1.5pt,shorten >=1.5pt,<->] (0.,0.)-- (17.5,0.);
\draw[shorten <=1.5pt,shorten >=1.5pt,<->] (0.,0.)-- (-35.,0.);
\draw(-8.75,8.505358310755952) circle (27.593543085191964);
\draw(-8.75,8.505358310755952) circle (12.202607098253445);
\draw[shorten <=1.5pt,shorten >=1.5pt,<->] (-20.263162558867325,12.548959167349107)-- (-35.,0.);
\draw[shorten <=1.5pt,shorten >=1.5pt,<->] (-20.263162558867325,12.548959167349107)-- (3.829443254744403,33.0646992505086);
\begin{small}
\draw [fill=white] (0.,0.) circle (1.5pt);
\draw[color=black] (0.12499708956780564,-2.5) node {$B$};
\draw [fill=white] (17.5,0.) circle (1.5pt);
\draw[color=black] (19.9,-1) node {$C$};
\draw[color=black] (8.571900784571307,-2.) node {$w$};
\draw [fill=white] (-35.,0.) circle (1.5pt);
\draw[color=black] (-37.5,-1) node {$A$};
\draw[color=black] (-20,-2.) node {$v$};
\draw [color=black] (-8.75,8.505358310755952) circle (1.5pt);
\draw[color=black] (-8.9,10.8) node {$Z$};
\draw[color=black] (-28.7,30.5) node {$c_1$};
\draw [fill=white] (3.829443254744403,33.0646992505086) circle (1.5pt);
\draw[color=black] (6.1,35.1) node {$E$};
\draw[color=black] (4.8,13.695096533418853) node {$c_2$};
\draw [fill=white] (-20.263162558867325,12.548959167349107) circle (1.5pt);
\draw[color=black] (-22.5,13.75) node {$D$};
\draw[color=black] (-29,8) node {$x_1$};
\draw[color=black] (-9,25) node {$x_2$};
\end{small}
\end{tikzpicture}\hfill
%\definecolor{uuuuuu}{rgb}{0.26666666666666666,0.26666666666666666,0.26666666666666666}
\begin{tikzpicture}[line cap=round,line join=round,>=latex,x=2.23,y=2.23]
\clip(-43.08416411948856,-15.544185487747239) rectangle (24.00374407352001,50.081758604203216);
\draw [shift={(-8.75,0.)}] (0,0) -- (90.:4.873213670194327) arc (90.:180.:4.873213670194327) -- cycle;
\draw[shorten <=1.5pt,shorten >=1.5pt,<->] (-8.75,0.)-- (17.5,0.);
\draw[shorten <=1.5pt,shorten >=1.5pt,<->] (-8.75,0.)-- (-35.,0.);
\draw(-8.75,16.781465191208603) circle (31.155739021306264);
\draw(-8.75,16.781465191208603) circle (16.781465191208603);
\draw[shorten <=1.5pt,shorten >=1.5pt,<->] (-35.,0.)-- (-22.889079190770428,7.742437942612188);
\draw[shorten <=1.5pt,shorten >=1.5pt,<->] (-22.889079190770428,7.742437942612188)-- (17.5,33.5629303824172);
\draw (-8.75,16.781465191208603)-- (-8.75,0.);
\fill(-10.77698966609729,2.0269896660972897) circle (0.16244045567314425);
\begin{small}
\draw [fill=white] (17.5,0.) circle (1.5pt);
\draw[color=black] (19.455411314671974,-2.9550501730785252) node {$C$};
\draw [fill=white] (-8.75,0.) circle (1.5pt);
\draw[color=black] (-8.646787516781984,-2.5) node {$B$};
\draw[color=black] (4.186008481396412,-2) node {$v$};
\draw [fill=white] (-35.,0.) circle (1.5pt);
\draw[color=black] (-37.5,-1) node {$A$};
\draw[color=black] (-24.078630805730686,-2) node {$v$};
\draw [fill=white] (-8.75,16.781465191208603) circle (1.5pt);
\draw[color=black] (-9.,19.5) node {$Z$};
\draw[color=black] (-31,41.63485490919969) node {$c_1$};
\draw[color=black] (-3.042591796058507,34.6) node {$c_2$};
\draw [fill=white] (17.5,33.5629303824172) circle (1.5pt);
\draw[color=black] (19.942732681691403,36.) node {$E$};
\draw [fill=white] (-22.889079190770428,7.742437942612188) circle (1.5pt);
\draw[color=black] (-21.8,10.9) node {$D$};
\draw[color=black] (-30.3,5.573073749761571) node {$x_1$};
\draw[color=black] (-2,23.766404785153775) node {$x_2$};
\draw[color=black] (-6.9,8) node {$r$};
\end{small}
\end{tikzpicture}

\parbox[t]{.48\linewidth}{Choose $c_1$ to be a sufficiently large circle through $A$ and $C$ and $c_2$
a concentric circle through $B$. Then $c_2$ cuts a chord $AE$ of length $u$ in
segments of lengths $x_1, x_2$. See also Exercise~\ref{sekantensatz}.} \hfill
\parbox[t]{.48\linewidth}{In particular, if we choose $v=w$ and $AE$ as diameter of $c_1$,
we have $r=\sqrt{(u/2)^2-v^2}$ and $|AZ|=u/2$. Hence we can read off the solutions of $x^2-ux+v^2=0$ 
and get the familiar expression $x_{1,2}=u/2\pm\sqrt{(u/2)^2-v^2}$.
} 
\end{center}
%This explains the
%three basic types of equations
%$$
%x^2+px=q, \quad x^2+q=px,\quad x^2=px+q.
%$$

%Griechen
%http://www.dm.uniba.it/~psiche/bas2/node4.html
%http://members.chello.at/gut.jutta.gerhard/geogl2.htm
% elemente Buch VU, Prop 28/29 (vd Waerden) vs: Sekantensatz
%http://aleph0.clarku.edu/~djoyce/elements/bookVI/propVI29.html

% The babylonian authors did not offer an explanation even in cases of two 
% positive solutions, only one is given. 
% is a cuneiform tablet AO 8862\footnote{Paris, Louvre, Antiquit\'es Orientales} from 
% Hammurabi dynasty (around 1700 BC). One of the exercise asks to solve (in modern
% notation)
% \begin{eqnarray*}
% x-y+xy&=&(3,3)_s = 3\times 60+3=183\\
% x+y&=&27
% \end{eqnarray*}
%which reduces to the quadratic equation $x^2-29x+210=0$. The babylonian
%%  [Ari_Ben-Menahem]_Historical_encyclopedia_of_natur(BookZZ.org).pdf Seite 2189 
%%  [Hans-Joachim_Waschkies__(auth.),__Jean_Christiani(BookZZ.org).pdf Seite 268
%%  Old Babylonian ?Algebra?, and What It Teaches Us about Possible Kinds of Mathematics Jens Høyrup Seite 21
%%  Neugebauer The exact sciences ... Seite 49
%%   https://fr.wikipedia.org/wiki/BM_13901ß

% http://experimentalmath.info/blog/2011/06/ancient-indian-square-roots/
% http://www-history.mcs.st-and.ac.uk/Projects/Pearce/Chapters/Ch6.html

The famous {\em Bakhshali manuscript}, the oldest passed down manuscript in
Indian mathematics, also contains examples and solutions of  quadratic
equations with numerical coefficients. This birch bark manuscript  is dated
between the 9th and the 12th century CE but is assumed to have origins 
going back to the 3rd or 4th century CE. 
% Brahmagupta
%%% Toledo Übersetzungsbewegung Cordoba Alfons der Weise von Castilien
In Chinese mathematics, the {\em Zhang Qiujian suanjing}, a
manual from the 5th century CE, not only discusses quadratic equations,
it also contains a consistent concept of negative numbers. See~\cite{ari} for 
more information.
%Liu Hui

An important milestone in the history of mathematics is
Al-Khwarizmi's algebra book {\em  Al-kit\={a}b al-mukhta\d{s}ar f\={\i} \d{h}is\={a}b 
al-\u{g}abr wa\hspace{-1pt}'l-muq\={a}bala\/} which is, in part, devoted to an
extensive discussion and classification of quadratic equations. 
Unlike in earlier work, Al-Khwarizmi, who lived between 780 and 850 CE, introduced the concept
of equation not only as a series of case by case problems to be solved, but
as an infinite class of problems with a common method of solution.
In particular, he identified the cases where the quadratic equation
has two, one or no solution.
He also demonstrated different methods of solution %, not algebraically,
%but geometrically, 
by means of the idea of equal area.
It is possible, that he was inspired by the translation
of Euclid's Elements by al-Hajjaj ibn Matar who, like  Al-Khwarizmi,
worked in the House of Wisdom in Baghdad (see~\cite[Chapter 11]{arab}).
Although most of the modern notation
was not yet available at that time and al-Khwarizmi still formulated
his problems and solutions in text form, his work marks the dawn of elementary algebra.
%https://books.google.ch/books?id=Dy897SeErOQC&pg=PA349&lpg=PA349&dq=Concise+Book+on+the+Computation+of+algebra+and+al-muqabala&source=bl&ots=aw_uSwIofP&sig=mGbitepEC2bQVaZyI72Gs4xLOzY&hl=de&sa=X&ved=0ahUKEwil5O29zcjRAhVMlxQKHS9GD4oQ6AEIKzAD#v=onepage&q=Concise%20Book%20on%20the%20Computation%20of%20algebra%20and%20al-muqabala&f=false

% griechen https://www.jstor.org/stable/2972623?seq=1#page_scan_tab_contents
% Muhammad al-Karkhi
% Abu Kamil Shuja (al-Hasib) 
% Liber embadorum
% Bar Hiyya
%Liu Hui Seite 740

So Descartes was, by far, not the first to solve quadratic equations.
However, instead of giving a verbal description of the algorithm
which yields the solution from the numerical values of the coefficients of 
a given single equation, he wrote down the solution as a condensed formula
involving the coefficients of the equation as variables. 
This progress can hardly be overestimated: 
The formula allows to analyze the solutions as {\em functions\/} of 
these coefficients, which in turn may depend on other
variables. Moreover, one may make use of the formula 
as an element in more complex computations.
% as part of more complex, comprehensive computations. 
Finally, %in particular,
a clear formal justification of the formula becomes possible, not
only a case by case verification of the numerical values at hand.

%The remarkable point is not, that Descartes found a way to compute
%the solutions of a quadratic equation: This algorithm has been known to
%the babylonians, arabs only for numbers and... The real innovation was
%that he found a way to write down the solutions in condensed form
%which could be analysed further and used to computer further computations.
%
%cannot be underestimated formula: Proof (before no proof only verification in each case)
%discussion (coefficients depend on functions) weiterrechnen.
\section{Friends of the quadratic formula}
Quadratic equations can not only be considered over $\mathbb C$ but also over a general field.
The way to compute the solutions of~(\ref{solabc}) from the coefficients
involves the four basic arithmetic operations and square roots.\vspace*{-3mm} 
\paragraph{Roots.}
In a general field $F$, $Q(F):=\{x^2\mid x\in F\}$ is the set of
squares. 
A function $\sqrt{\textcolor{white}p } :Q\to F$ with 
the property $(\sqrt x)^2=x$ is called a square root.
The polynomial $p(x)=x^2-r$ has at most two zeros. For $r\in Q(F)$
the two zeros are $\pm\sqrt r$, and these two values 
coincide if $r=0$ or if $F$ has characteristic 2.
If $F$ is algebraically closed all elements are squares.
By virtue of the Frobenius automorphism the same is true for every finite field of  
characteristic 2.
And in a finite field of odd characteristic
there are as many nonzero squares as
nonsquares. In that case, products of two squares or of two nonsquares are 
squares, a product of a square and a nonsquare is a nonsquare.
%http://math.stackexchange.com/questions/1938366/in-a-finite-field-product-of-non-square-elements-is-a-square

The axiom of choice, or more precisely the axiom of choice $C_2$ for families of
pair sets, ensures the existence of a square root
for a general field. Actually, $C_2$ is equivalent to the fact, that
in each field there exists a root on its set of squares (see%~\cite[Chapter 5, Related Results]{halbeisen} and
~\cite[Proposition 4.3]{hhllls}). 
However, there are cases where a root may be defined
due to additional structures:
If $F$ is an ordered field, e.g. $\mathbb R$, one may always choose the larger of two possible values for the root of a square.
In the case $F=\mathbb C$ one can use its  structure of a Riemann surface
to define a root by identifying an analytical principal branch of the root function.

It is worth noticing that
the solution formula~(\ref{solabc}) works with {\em any\/} choice of a root function
and in {\em any\/} field of  characteristic $\neq 2$. Let us have a
brief look at the special case of a field of characteristic $2$.
\paragraph{Quadratic equations in a field of characteristic 2.}
We consider~(\ref{abc}) in a field $F$ of characteristic 2. If $b=0$, 
then~(\ref{abc}) is equivalent to $x^2=c/a$. And if $c/a\in Q(F)$, 
or equivalently $ac\in Q(F)$,
the unique solution in $F$ is $\sqrt{c/a}$. If $ac\notin Q(F)$, there is no solution in $F$.

If $b\neq 0$, we rewrite~(\ref{abc}) in the form 
$$
\bigl(\frac{ax}b\bigr)^2+\frac{ax}b+\frac{ac}{b^2}=0
$$
which, for $y:=\frac{ax}b$, takes the form 
\begin{equation}
y^2+y=\frac{ac}{b^2}.\label{artin}
\end{equation}
If we consider $F$ as  a vector space over its prime field $\mathbb F_2$,
then the Artin-Schreier map $F\to F, y\mapsto y^2+y$,
is linear. Its kernel $K$ consists of  the elements $z\in F$ for which $z^2=z$,
i.e., $K=\mathbb F_2$. Hence,~(\ref{artin}) has a solution $y$ iff $ac/b^2$
is an element of the image of the Artin-Schreier map, and then $y+1$
is the second solution. In terms of the original variable $x$, we get:
(\ref{abc}) with $b\neq 0$ has no solution if $ac/b^2$ does not lie in the image
of the Artin-Schreier map, and it has the two solutions $\frac{by}a$ and 
$\frac{b(y+1)}a$ if $y$ is a preimage of $ac/b^2$ under the Artin-Schreier map.
For a detailed discussion, see~\cite{char2}.
%https://www.staff.uni-mainz.de/pommeren/MathMisc/QuGlChar2.pdf
%
%
\paragraph{Characteristic \boldmath$\neq 2$.} In a general field $F$ of
characteristic $\neq 2$,
the number of solutions of~(\ref{abc}) is governed by the
discriminant $D:=b^2-4ac$:
For coefficients $a,b,c\in F$, equation (\ref{abc}) has
no solution in $F$
if the discriminant $D$ is a nonsquare.
If $D$ is a square, 
(\ref{abc}) has one solution if $D=0$ %or if $F$ has characteristic $2$,
otherwise two solutions. 

But let us now concentrate to the case $F=\mathbb C$
and discuss alternative ways to represent the solutions of~(\ref{abc}).

When starting from the normalized form
\begin{equation}\label{normal}
x^2-2px+q=0, \quad p,q\in\mathbb C,
\end{equation}
the quadratic formula takes the form
\begin{equation}\label{solnormal}
x_{1,2}=p\pm\sqrt{p^2-q}.
\end{equation}
Lesser-known is the formula
\begin{equation}\label{muller}
x_{1,2}=\frac{-2c}{b\pm\sqrt{b^2-4ac}}
\end{equation}
for the solutions of~(\ref{abc}). 
%This form is obtained by
%dividing~(\ref{abc}) by $x^2$, solving the resulting equation for $1/x$
%by formula~(\ref{solabc}), and then taking the reciprocal.
This variant follows from~(\ref{solabc}) by observing, that
$x_1x_2=c/a$.
The formula (\ref{muller}) yields only one solution if $c=0$,
but otherwise it
%has the advantage of being 
%particularly robust: It 
gives a correct answer
even in the degenerate case  $a=0$ (i.e.~a linear equation), provided $b\neq 0$. 
For this reason~(\ref{muller}) is used in Muller's iterative method~\cite{muller}
for numerically finding roots of a function $f$ by approximating $f$
with  parabolas interpolating the graph  in three points.

For real coefficients, trigonometric and hyperbolic methods 
for solving quadratic equations
are also circulating: In contrast to~(\ref{solabc}), (\ref{solnormal}) and~(\ref{muller})
these formulas are numerically stable (see~\cite{fortran}, \cite{pitfalls}, and  
the discussion in Section~\ref{newfriend}).
To explain the idea we assume, without loss of generality, that
$a>0$ in~(\ref{abc}). Moreover, we sort out the cases $b=0$ 
(with the solutions $x_{1,2}=\pm \sqrt{-c/a}$)
and $c=0$ (with the solutions $x_1=0,x_2=-b/a$).
Otherwise, we distinguish the cases $c>0$ and $c<0$:

{\em Case 1.} If $c<0$ we bring equation~(\ref{abc})
by the scaling $x=y\sqrt{-\frac{c}{a}}$ to the normalized form
\begin{equation}\label{trigo}
y^2+dy-1=0, \text{ with $d=\frac{b}{\sqrt{-ac}}$}.
\end{equation}
Using the double-angle formula
$$
\tan 2\alpha=\frac{2y}{1-y^2}\quad\text{for $y=\tan\alpha$}
$$
we see that $y_1=\tan\alpha$  is a solution of~(\ref{trigo}) if 
$\alpha$ is chosen such that $\tan2\alpha=\frac{2\sqrt{-ac}}b$. According to
Vi\"eta's Theorem, the product of the solutions of~(\ref{trigo}) is $y_1y_2=-1$,
i.e.,~$y_2=-\cot\alpha$. Summarizing, we get as solutions of~(\ref{abc})
$$
x_1=\sqrt{-\frac{c}{a}}\tan\alpha,\quad 
x_2=-\sqrt{-\frac{c}{a}}\cot\alpha, \quad\text{with $\tan2\alpha=\frac{2\sqrt{-ac}}b$.}
$$
{\em Case 2.} If $c>0$ we first have to sort out the case
$b^2-4ac= 0$ with double root $x_1=x_2=-\frac{b}{2a}$.
Otherwise, we bring equation~(\ref{abc})
by the scaling $x=y\sqrt{\frac{c}{a}}$\footnote{Alternatively, one
may use $x=iy\sqrt{\frac{c}{a}}$ to arrive again at the form $y^2+dy-1=0$ like in the first case.} to the normalized form
\begin{equation}\label{hyper}
y^2+dy+1=0, \text{ with $d=\frac{b}{\sqrt{ac}}$}.
\end{equation}
Recall the double-angle formula for the hyperbolic tangent
\begin{equation}\label{hyper1}
%\tan \alpha=\frac{2y}{1-y^2}\text{ for $y=\tan\frac\alpha2$}
%\text{ and } 
\tanh2\alpha=\frac{2y}{1+y^2} \quad\text{for }y=\tanh\alpha.
\end{equation}
%which is valid outside the poles of $\mathbb C\to\mathbb C, \alpha\mapsto \tanh2\alpha$,
%i.e., for $\alpha\in\mathbb C\setminus\{\frac{i\pi}4+\frac{k\pi i}2:k\in\mathbb Z\}$.
The range of the function $\tanh$ is $\mathbb C\setminus\{-1,1\}$.
%Since $b^2-4ac\neq 0$, we can choose $\alpha$ such that
%$\tanh2\alpha=-\frac{2\sqrt{ac}}b\neq \pm 1$: 
Since $b^2-4ac\neq 0$, we have $\frac{2\sqrt{ac}}b\neq \pm1$
and we can choose $\alpha$ such that $\tanh2\alpha=-\frac{2\sqrt{ac}}b$,
and actually $\alpha\in\mathbb R$ iff $b^2-4ac> 0$.
%
%If $b^2-4ac< 0$ we have $\left|\frac{2\sqrt{ac}}b\right|>1$, and there exist
%$\alpha\in\mathbb C\setminus\mathbb R$ with $\tanh2\alpha=-\frac{2\sqrt{ac}}b$. 
Then, by~(\ref{hyper1}) it follows 
that $y_1=\tanh\alpha$  is a solution of~(\ref{hyper}).  Again by Vi\"eta's Theorem, 
the product of the solutions of~(\ref{hyper}) is $y_1y_2=1$,
i.e.,~$y_2=\coth\alpha$. Putting all together, we find as solutions of~(\ref{abc})
in the considered case
$$
x_1=\sqrt{\frac{c}{a}}\tanh\alpha,\quad x_2=\sqrt{\frac{c}{a}}\coth\alpha, \quad\text{with $\tanh2\alpha=-\frac{2\sqrt{ac}}b$.}
$$
%Notice that  $b^2-4ac< 0$ corresponds to $\left|\frac{2\sqrt{ac}}b\right|>1$ 
%leading to $\alpha\in\mathbb C\setminus\mathbb R$ and ultimately to
%two complex conjugate solutions.

It may look strange at first sight to employ trigonometric or hyperbolic functions
to solve a quadratic equation. However, recall, that the use of trigonometric functions
becomes inevitable  for the cubic and quartic equation in the {\em casus irreducibilis}.
And for the general quintic equation, even trigonometric functions
are no longer strong enough: Hermite used the theory of elliptic modular functions
and the Jacobi theta function to solve the quintic 
(see~\cite{hermite}). Later, Klein related the symmetries of the icosahedron 
with Galois theory to show systematically, why elliptic modular functions
appear naturally in the solution of the quintic. He also devised
a solution using hypergeometric functions
(see~\cite{klein}).

\section{Meet a new friend}\label{newfriend}
Here, we want to propose an alternative quadratic formula 
which stands out by its symmetry. 
\begin{prop}
Let the quadratic equation be given by the homogeneous normal form
\begin{equation}\label{alt}
x^2-4ux+4v^2=0,\quad u,v\in\mathbb C.
\end{equation}
%It refers to the normal form
%\begin{equation}\label{alt}
%x^2-4ux+4v^2=0,\quad u,v\in\mathbb C.
%\end{equation}
Then its solutions can be expressed by
\begin{equation}\label{solalt}
x_{1,2}=(\sqrt{u+v}\pm\sqrt{u-v})^2
%\end{equation}
\text{ or equivalently }
x_{1,2}=(\sqrt{u-v}\pm\sqrt{u+v})^2.
\end{equation}
For real coefficients $u,v^2\in\mathbb R$, (\ref{alt}) has\vspace*{-1mm}
\begin{itemize}\itemsep=1mm
\item[(a)] a  double root (namely $2u$) iff $u^2=v^2$,
\item[(b)] two positive distinct real roots iff  $u^2>v^2$,
\item[(c)] two real roots of opposite sign iff $v^2<0$.
\end{itemize}
\end{prop}
{\em Remark.} It is easy to convert the standard form (\ref{abc}) 
to the homogeneous form (\ref{alt}), however
it requires a root to determine $v$. But for dimensional reasons,
problems originating from physics or geometry often present themselves
already in the form (\ref{alt}). Actually, Descartes' first example
at the very hour of birth of the quadratic formula, was geometric and of the homogeneous form (\ref{alt})
with coefficients $u=\frac a4, v=\frac b2$: The following facsimile shows the 
corresponding excerpt of page 303 in~\cite{descartes}.
\begin{center}
\includegraphics[width=.7\linewidth]{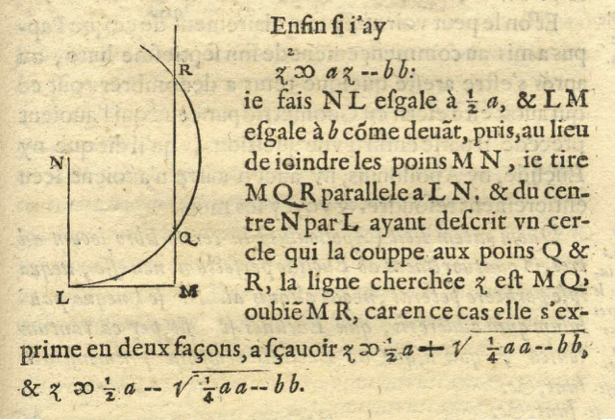}
\end{center}
%Here, we want to propose an alternative quadratic formula which stands out by
%its symmetry. It refers to the normal form~(\ref{alt})
%and reads as follows:
%Then, the solution of~(\ref{alt}) can be written in the following form
%which stands out by its symmetry:
{\em Proof of the Proposition.}
The classical formula~(\ref{solabc}) gives initially 
\begin{equation}\label{initially}
x_{1,2}=2u\pm2\sqrt{u^2-v^2}.
\end{equation}
This establishes (a)--(c) in the case of real coefficients $u,v^2$.
In general, if $u,v^2\in\mathbb C$, we use that for any root function on $\mathbb C$ the product
$\sqrt{u+v}\sqrt{u-v}$ is either plus or minus $\sqrt{(u+v)(u-v)}$ and get by~(\ref{initially})
\begin{align}
x_{1}&=2u+2\sqrt{u^2-v^2}=u+v+u-v+2\sqrt{(u+v)(u-v)}=\notag\\
&=u+v+u-v\pm2\sqrt{u+v}\sqrt{u-v}=(\sqrt{u+v}\pm\sqrt{u-v})^2\notag%\tag*{$\Box$}
\end{align}
with either the plus or the minus sign. The result for $x_2$ is the
same expression with the opposite sign.

Alternatively, the correctness of~(\ref{solalt}) can readily be 
checked, without transforming the product $\sqrt{u+v}\sqrt{u-v}$,
by plugging it into~(\ref{alt}). \hfill$\Box$
%Now, for the solutions of the quadratic equation in normal form~(\ref{alt}) 
%we propose the following formula 
%which stands out by its symmetry:
%\begin{equation}\label{solalt}
%x_{1,2}=(\sqrt{u+v}\pm\sqrt{u-v})^2
%\end{equation}
%or, equivalently,
%$$x_{1,2}=(\sqrt{u-v}\pm\sqrt{u+v})^2.$$
% The correctness of the formula is readily checked by plugging it into~(\ref{alt}).
% (Observe, that $\sqrt r\sqrt s\neq\sqrt{rs}$ in general.)
%expanding the 
%square and comparing to~(\ref{solnormal}). 
% By looking at~(\ref{solalt}),
%we read off that, for real coefficients $u,v^2\in\mathbb R$, (\ref{alt}) has
% \begin{itemize}
% \item a  double root (namely $2u$) iff $u^2=v^2$,
% \item two distinct real roots iff  $u^2>v^2$,
% \item two real roots of opposite sign iff $v^2<0$.
% \end{itemize}

A benefit of the new formula~(\ref{solalt}) compared to
the classical version~(\ref{solabc}) is its numerical stability. To illustrate
this property, we compare the classical formula
with the new formula  when applied to the model problem
in~\cite{pitfalls}:
$$
x^2-10^nx+1=0.
$$
For simplicity we take $n>1$ even.
Then % , as explained in~\cite{pitfalls}, 
the solution obtained by~(\ref{solabc})
rounded to $2n$ significant decimals, are 
\begin{align*}
x_1&=10^n-10^{-n}&&\hspace*{-13mm}\text{with a relative error of the order $10^{-4n}$}\\
x_2&=10^{-n}&&\hspace*{-13mm}\text{with a relative error of the order $10^{-2n}$.}
\end{align*}
However, if we repeat the calculation rounded to $2n-1$ significant decimals, 
we are victims of catastrophic cancellation and get
\begin{align*}
x_1&=10^n&&\hspace*{-13mm}\text{with a relative error of the order $10^{-2n}$}\\
x_2&=0&&\hspace*{-13mm}\text{with a relative error of $100\%$.}
\end{align*}
If we use~(\ref{solalt}) instead, the situation is considerably better, as it
still yields satisfactory solutions if the calculations are carried out 
rounded to only $n$ significant decimals:
\begin{align*}
x_1&=10^n&&\hspace*{-13mm}\text{with a relative error of the order $10^{-2n}$}\\
x_2&=10^{-n}&&\hspace*{-13mm}\text{with a relative error of the order $10^{-2n}$.}
\end{align*}

The following example taken from a physics textbook~\cite[Section 14.3, Problem 8]{budde}
shows how nicely~(\ref{solalt}) may turn out compared to the classical
formulas~(\ref{solabc}) or~(\ref{solnormal}):

{\bf Exercise.} A stone is dropped from rest into a well. $t$ seconds later, the sound
of the splash is heard. Ignoring air resistance, how far has the stone fallen before 
hitting the water?

{\bf Solution.} Let $x$ denote the depth of the well, $c$ the velocity of sound,
 and $g$ the acceleration of fall (we neglect air resistance).
The time $t$ is the sum of the time $t_{1}=\sqrt{\frac{2x}{g}}$ for the 
stone to fall and the travel time $t_2=\frac{x}{c}$ of the sound. Squaring
$t_1=t-t_2$ and rearranging terms leads to the quadratic equation
\begin{equation}\label{well}
x^2-2x\frac{c(c+gt)}{g}+t^2c^2=0.
\end{equation}
The alternative formula~(\ref{solalt}) yields the particularly compact
expression
\begin{equation}\label{stone}
x_1=\frac{c^2}{2g}\Bigl(
\sqrt{1+\frac{2gt}{c}}-1
\Bigr)^2.
\end{equation}
Compare this to the rather messy expression you get from the classical formula~(\ref{solnormal})!
Let us briefly discuss the result: A falling stone reaches after time $\tau=\frac cg$
the velocity of sound. During that time $\tau$, it travelled the distance $\frac{c^2}{2g}$
which is precisely the first factor in~(\ref{stone}). The second factor, the squared bracket, is a
dimensionless number. The Taylor expansion in powers of $1/c$ of~(\ref{stone}) is
$$
x_1=\frac{gt^2}2-\frac{g^2t^3}{2c}+O(1/c^2).
$$
Here, the first term is the answer to the exercise in the case of infinite propagation of sound,
the second term gives the first order correction (in $1/c$) which yields reasonably
good results in practice.

Observe, that the second solution $x_2$ does not solve the
given physical problem. However,
the reader is invited to figure out what physical
interpretation the second solution of~(\ref{well}) has.

\section*{Exercises}
\begin{enumerate}
\item Solve Descartes' problem $z^2=az-b^2$ by the alternative formula~(\ref{solalt}).
Interpret the result geometrically.

\item Let $F$ be a field of characteristic $\neq 2$. Prove: The solutions of the 
quadratic equation $x^2-px+q=0$ are squares in $F$ if and only if $q=r^2$ is a square,
and $p\pm 2r$ are squares.

\item Let $a,b\in\mathbb Q$ be positive such that $a^2>b$, $b\notin Q(\mathbb Q)$.
Show that $x,y\in\mathbb Q$ with 
$$
\sqrt{a\pm\sqrt b}=\sqrt x\pm\sqrt y
$$
exist iff $a^2-b\in Q(\mathbb Q)$. Find $x$ and $y$ in this case.
%% Dietmar Herrmann: Die antike Mathematik, Seite 152, Seite 261

\item\label{sekantensatz} Solve the quadratic equation $x^2-ux-vw=0$ with the 
Intersecting Secants Theorem (Theorem 36 in Book III of Euclid's {\em Elements}).
% x -> -x to change the sign of u

\item A company sells olives in large cardboard boxes folded from squares of sidelength $a$:
\begin{center}
\begin{tikzpicture}[line cap=round,line join=round,>=latex,x=0.5cm,y=0.5cm]
\clip(-1.5062229967591512,-3.211564822190765) rectangle (18.118308564442163,5.153568235371656);
\fill[fill=black,fill opacity=0.1] (0.,0.) -- (2.,0.) -- (2.,-2.) -- (5.,-2.) -- (5.,0.) -- (7.,0.) -- (7.,3.) -- (5.,3.) -- (5.,5.) -- (2.,5.) -- (2.,3.) -- (0.,3.) -- cycle;
\fill[fill=black,pattern=north east lines,pattern color=black] (2.,5.) -- (2.,3.) -- (0.,3.) -- (0.,5.) -- cycle;
\fill[fill=black,pattern=north east lines,pattern color=black] (5.,5.) -- (5.,3.) -- (7.,3.) -- (7.,5.) -- cycle;
\fill[fill=black,pattern=north east lines,pattern color=black] (5.,0.) -- (7.,0.) -- (7.,-2.) -- (5.,-2.) -- cycle;
\fill[fill=black,pattern=north east lines,pattern color=black] (2.,-2.) -- (2.,0.) -- (0.,0.) -- (0.,-2.) -- cycle;
\fill[fill=black,fill opacity=0.1] (13.,0.) -- (17.,0.) -- (17.,2.) -- (13.,2.) -- cycle;
\fill[fill=black,fill opacity=0.1] (17.,0.) -- (18.,1.) -- (18.,3.) -- (17.,2.) -- cycle;
\fill[fill=black,fill opacity=0.25] (13.,2.) -- (14.,3.) -- (18.,3.) -- (17.,2.) -- cycle;
\fill[fill=black,fill opacity=0.25] (14.,3.) -- (13.,2.) -- (14.,2.) -- cycle;
\draw (0.,0.)-- (2.,0.);
\draw (2.,0.)-- (2.,-2.);
\draw (2.,-2.)-- (5.,-2.);
\draw (5.,-2.)-- (5.,0.);
\draw (5.,0.)-- (7.,0.);
\draw (7.,0.)-- (7.,3.);
\draw (7.,3.)-- (5.,3.);
\draw (5.,3.)-- (5.,5.);
\draw (5.,5.)-- (2.,5.);
\draw (2.,5.)-- (2.,3.);
\draw (2.,3.)-- (0.,3.);
\draw (0.,3.)-- (0.,0.);
\draw (0.,5.)-- (0.,-2.);
\draw (0.,-2.)-- (7.,-2.);
\draw (7.,-2.)-- (7.,5.);
\draw (7.,5.)-- (0.,5.);
\draw (2.,5.)-- (2.,3.);
\draw (2.,3.)-- (0.,3.);
\draw (0.,3.)-- (0.,5.);
\draw (0.,5.)-- (2.,5.);
\draw (5.,5.)-- (5.,3.);
\draw (5.,3.)-- (7.,3.);
\draw (7.,3.)-- (7.,5.);
\draw (7.,5.)-- (5.,5.);
\draw (5.,0.)-- (7.,0.);
\draw (7.,0.)-- (7.,-2.);
\draw (7.,-2.)-- (5.,-2.);
\draw (5.,-2.)-- (5.,0.);
\draw (2.,-2.)-- (2.,0.);
\draw (2.,0.)-- (0.,0.);
\draw (0.,0.)-- (0.,-2.);
\draw (0.,-2.)-- (2.,-2.);
\draw [shorten >=-4pt,-{[scale=1.2]>},decorate,decoration=snake] (8.517971510584903,1.006297626664551) -- (11.517971510584903,1.006297626664551);
%\draw [->,decorate,decoration=snake] (8.517971510584903,1.006297626664551) -- (11.517971510584903,1.006297626664551);
\draw (13.,0.)-- (17.,0.);
\draw (17.,0.)-- (17.,2.);
\draw (17.,2.)-- (13.,2.);
\draw (13.,2.)-- (13.,0.);
\draw (17.,0.)-- (18.,1.);
\draw (18.,1.)-- (18.,3.);
\draw (18.,3.)-- (17.,2.);
\draw (17.,2.)-- (17.,0.);
\draw (13.,2.)-- (14.,3.);
\draw (14.,3.)-- (18.,3.);
\draw (18.,3.)-- (17.,2.);
\draw (17.,2.)-- (13.,2.);
\draw (14.,3.)-- (13.,2.);
\draw (13.,2.)-- (14.,2.);
\draw (14.,2.)-- (14.,3.);
\draw [->] (-1.,5.) -- (-1.,-2.);
\draw [->] (-1.,-2.) -- (-1.,5.);
\draw [->] (0.,-3.) -- (7.,-3.);
\draw [->] (7.,-3.) -- (0.,-3.);
\draw [->] (8.,0.) -- (8.,-2.);
\draw [->] (8.,-2.) -- (8.,0.);
\draw [dash pattern=on 3pt off 3pt] (2.,3.)-- (5.,3.);
\draw [dash pattern=on 3pt off 3pt] (5.,3.)-- (5.,0.);
\draw [dash pattern=on 3pt off 3pt] (5.,0.)-- (2.,0.);
\draw [dash pattern=on 3pt off 3pt] (2.,0.)-- (2.,3.);
\begin{small}
\draw[color=black] (-1.3,1.5) node {$a$};
\draw[color=black] (3.5,-2.7) node {$a$};
\draw[color=black] (8.33,-1.) node {$x$};
\end{small}
\end{tikzpicture}
\end{center}
The selling price for olives is $s$ Dollar per unit volume. The price for one
box is therefore $p(x)=sx(a-2x)^2$.
The cost for the cardboard is $c$ Dollar per unit area. The
cost per box is therefore $q(x)=c(a^2-4x^2)$. Optimize the profit
$p(x)-q(x)$ per box. Use the alternative formula~(\ref{solalt}).
\end{enumerate}

{\bf Acknowledgement.} The author would like to thank Hans Peter Dreyer for pointing out
to him the nice exercise of the falling stone, and 	
Jacques G\'elinas for the remark about the numerical stability.
%The author is also grateful for the critical remarks and hints of the referees
%which helped to improve the first version of the article.
%\bibliography{quadratic_formula_1}
%\bibliographystyle{plain}

{\small Norbert Hungerb\"uhler\\
Department of Mathematics\\
ETH Z\"urich\\
8092 Z\"urich, Switzerland\\
norbert.hungerbuehler@math.ethz.ch
}\\[10mm]

 \end{document}